%
%
%
%


\magnification=\magstep1                        
\input amstex
\documentstyle{amsppt}
\topmatter
\nopagenumbers
\vsize=8.605truein
\hsize=6.55truein
\voffset0.3truein
\title
quasiconformal homeomorphisms on cr 3-manifolds with symmetries
\endtitle
\author
PUQI TANG
\endauthor
\abstract
An extremal quasiconformal homeomorphisms in a class of homeomorphisms
between two CR 3-manifolds is an one which has the least conformal distortion
among this class. This paper studies extremal quasiconformal homeomorphisms
between CR 3-manifolds which admit transversal CR circle actions. Equivariant
$K$-quasiconformal homeomorphisms are characterized by an area-preserving
property and the $K$-quasiconformality of their quotient maps on the
spaces of $S^1$-orbits. A large family of invariant CR structures on $S^3$
is constructed so that the extremal quasiconformal homeomorphisms among
the equivariant mappings between them and the standard structure are
completely determined. These homeomorphisms also serve as examples showing
that the extremal quasiconformal homeomorphisms between two invariant
CR manifolds are not necessarily equivariant. 
\endabstract
\keywords 
contact and CR structures, quasiconformality, variation, $S^1$-action 
\endkeywords
\subjclass
32G07, 30C70
\endsubjclass 
\address
Department of Mathematics, Purdue University, West Lafayette, IN 47907, USA
\endaddress
\email
tang\@math.purdue.edu \newline \newline
\endemail
\endtopmatter
\document

\heading
1. Introduction
\endheading

Given an oriented, compact, smooth surface $R$ of genus $>1$, divide all complex
structures on $R$ into equivalence classes so that two structures are in
the same class if and only if there is a conformal homeomorphism between
them which is homotopic to the identity. Teichm\"uller's theorem says that
for any two complex structures $S_1$ and $S_2$ on $R$, among all quasiconformal
homeomorphisms homotopic to the identity, there is an unique homeomorphism
which minimizes the conformal distortion with respect to $S_1$ and $S_2$, and
this extremal quasiconformal homeomorphism can be characterized in terms of
certain holomorphic quadratic differentials [2]. The maximal
dilatation of extremal quasiconformal homeomorphism measures how
different the class $[S_1]$ is from the class $[S_2]$. Since these fundamental
results have been established, Teichm\"uller space, the space of all equivalence
classes, became one of the most important objects of research in complex
analysis. Comprehensive literatures on Teichm\"uller theory include Abikoff's
[1], Zhong Li's [14] and Nag's [17].

Lempert proposed an analogous problem in the setting of Cauchy-Riemann (CR)
manifolds as follows [13]. Given two CR structures on a 3-dimensional contact
manifold, describe the quasiconformal homeomorphisms that have the least
conformal distortion with respect to these two CR structures. These
homeomorphisms, if exist, are said extremal. Their maximal dilatation
measures the nonisomorphism of the two CR structures. A Teichm\"uller type
distance between the two CR manifolds is defined by the infimum of the
logarithms of the maximal dilatations of all quasiconformal
homeomorphisms between them. This can be regarded as a variational approach
to the embeddability of an abstract CR structure. If the distance between
an abstract CR structure and an embeddable CR structure is zero and is also
realized, then the abstract CR structure is conformally equivalent to the
embedded one. We were able to prove that conformal equivalence implies CR
equivalence for embeddable CR structures, and we conjecture this holds for
general CR structures. Otherwise, one would like to know how far this CR
structure is from the space of all embeddable structures.

The concept of quasiconformality is classically given on Riemann surfaces
and Riemannian manifolds. It is a major machinery applied in Teichm\"uller
theory. Mostow introduced it for symmetric spaces of real rank one, which
include the Heisenberg groups [16]. Later Kor\'anyi and Reimann generalized
notion of quasiconformality to strongly pseudoconvex CR manifolds [10].

We will study extremal quasiconformal homeomorphisms between smooth, compact,
strongly pseudoconvex CR manifolds of dimension 3. In this paper, we shall
mostly consider CR manifolds that admit a transversal CR action of $S^1$,
in particular, the 3-sphere $S^3$ with the standard circle action. We remark
that these CR structures are always embeddable ([6] [12]); if the underlying
contact manifold is $S^3$, they can even be embedded into ${\Bbb C}^2$ as
circular hypersurfaces [6].

There are two basic questions here. The first question is
whether an extremal quasiconformal homeomorphism between two
$S^1$-invariant CR structures is $S^1$-equivariant. The second question
is what is the characterization of equivariant quasiconformal homeomorphisms.

The space of $S^1$-orbits of an invariant CR manifold is a surface
with a complex structure induced from the CR structure. An equivariant
homeomorphism between two $S^1$-invariant CR manifolds defines a
quotient homeomorphism between the corresponding Riemann surfaces.
In this paper we prove that an equivariant $K$-quasiconformal
homeomorphism is characterized by an area-preserving property and
$K$-quasiconformality of its quotient homeomorphism (Theorem 3.5, 3.6).
This answers the second question. We also develop the first and
second variation of the conformal distortion on $S^3$ (Proposition 5.1, 5.3).
The method to compute the variation on $S^3$ works on any CR 3-manifolds.
Then we construct a family of smooth $S^1$-invariant CR structures on
$S^3$ so that no extremal quasiconformal homeomorphism between these CR
structures and the standard CR structure is $S^1$-equivariant
(Theorem 6.1). Thus we show that circular symmetry is broken for
extremal quasiconformal homeomorphisms between these $S^1$-invariant CR
structures.

Recently we found that in certain situations an extremal quasiconformal
homeomorphism in a homotopy class must be equivariant. There
the extremal homeomorphisms have behavior analogous to Teichm\"uller 
transformations on Riemann surfaces. Details will appear in a forthcoming paper. 

\subhead
Acknowledgements
\endsubhead
This work was done under the guidance of L\'aszl\'o Lempert, my major
professor. I would like to express my deep gratitude to him for
leading me to this area and for his numerous invaluable suggestions, pleasant
teaching and persistent patience. I would like to thank Zhong Li for
teaching me Teichm\"uller theory and sending me his interesting book
on this subject. I also thank Kor\'anyi, Gehring, Mostow and Reimann
for providing me related papers some of which are main references of
this paper. 

\heading
2. Quasiconformal Homeomorphisms and Contact Flows 
\endheading
 
Let $M$ be a 3-dimensional, connected, smooth, contact manifold with a 
smooth non-degenerate contact form $\eta$. Denote the contact bundle by
$HM \triangleq \text{Ker }\eta$. Let $J_0: HM \to HM$ be a smooth endomorphism
such that $J_0^2=-\text {id}$. Thus $J_0$ is a smooth complex structure
on $HM$ which defines a strongly pseudoconvex CR structure on $M$. The
corresponding CR manifold is denoted by $M_0$.

Call the orientation of $M$ given by $d\eta\wedge\eta\ne 0$ positive 
and the orientation of $HM$ given by $d\eta\vert_{HM}$ positive. Note if
$\eta^\prime =\lambda\eta$ with a function $\lambda\ne 0$ is another
contact form, the orientation of $M$ given by $d\eta^\prime\wedge\eta^\prime
=\lambda^2 d\eta\wedge\eta$ is positive. The orientation of $HM$ given by
$d\eta^\prime\vert_{HM}=\lambda\, d\eta\vert_{HM}$ is either positive when
$\lambda >0$ or negative when $\lambda <0$.

Let $X\ne 0$ be a local section of $HM$, then $X$ and $J_0X$ are linearly
independent. $d\eta$ is nondegenerate on $HM$, so $\langle d\eta,
X\wedge J_0X\rangle \ne 0$. We say the CR structure of $M_0$ is positively
(or negatively) oriented with respect to $\eta$ if $\langle d\eta,X\wedge
J_0X\rangle >0$ (or $<0$). Note 
$$
\langle d\eta\wedge\eta, X\wedge J_0X\wedge [J_0X,X]\rangle =(\langle d\eta,
X\wedge J_0X\rangle)^2 >0.
\tag2.1
$$
Hence $X, J_0X, [J_0X, X]$ is always a positively oriented frame no matter the
CR structure is positively oriented or not. 

A differentiable curve on $M$ is called Legendrian if its tangent vector
at each point is in the contact bundle $HM$. Let $U\i M$ be an open set,
$\Gamma$ be a contact fibration of $U$, i.e., $\Gamma$ is a smooth fibration
of $U$ consisting of smooth Legendrian curves. A subfamily $\Gamma_1$ of a
contact fibration $\Gamma$ of $U$ is said to be of measure zero if for any
smooth surface $S$ which is transversal to each $\gamma\in\Gamma$ and any
smooth area form $\omega$ on $S$
$$
\displaystyle \int_{\{ S\cap\gamma |\gamma\in\Gamma_1\} }\omega =0. \tag2.2
$$

Assume that $M_1$ is another smooth, strongly pseudoconvex CR manifold with
the same underlying contact manifold $M$ and a complex structure $J_1$ on
$HM$. A homeomorphism $f: M_1\to M_0$ is said to be ACL (absolutely continuous
on lines) if for any open set $U\i M$ and contact fibration $\Gamma$ of $U$,
$f$ is absolutely continuous along all curves in $\Gamma$ except for a
subfamily of $\Gamma$ of measure zero.

For $j=0,1$, let $HM_j$ denote $HM$ endowed with the CR structure $J_j$.
Take any Hermitian metric on $HM_j$ with respect to $J_j$. Denote by
$|\cdot |_j$ the corresponding norm on $HM_j$.

\proclaim {Definition 2.1} (i)  A homeomorphism $f:M_1 \to M_0$ is
$K$-quasiconformal if
\roster
\item $f$ is ACL;

\item $f$ is differentiable almost everywhere and its differential
$f_*$ preserves the contact bundle; and

\item the maximal dilatation $\dsize K= K(f)= \underset{q\in M_1}
\to{\text{\rm ess sup}}\, K(f)(q) <\infty $, where
$$
K(f)(q)=\frac {\displaystyle \max _{X \in H_qM_1,|X|_1=1}
|f_*X|_0}{\displaystyle \min_{X \in H_qM_1,|X|_1=1}|f_*X|_0}.
\tag2.3
$$
is the dilatation of $f$ at $q \in M_1$.
\endroster
(ii) A 1-quasiconformal homeomorphism $f:M_1 \to M_0$ is called conformal.
If such a conformal homeomorphism exists, $M_1$ and $M_0$ are said conformally
equivalent.
\endproclaim

\demo{Remark} (1) For any $q\in M$, $j= 0,1$, $\dim_{\Bbb C}H_qM_j=1$, so
any two Hermitian metric on $H_qM_j$ are scalar multiples of each other.
Hence the value of $K(f)(q)$ is independent of the choices of the
Hermitian metrics.

(2) A $C^1$ homeomorphism is conformal if and only if it is CR. When both
$M_0$ and $M_1$ are smooth and embeddable into ${\Bbb C}^2$, a homeomorphism
$f: M_1\to M_0$ is conformal if and only if it is smooth and CR. A proof
to this will be given in a forthcoming paper.

(3) On the standard 3-sphere, Kor\'anyi and Reimann gave an analytic
definition of quasiconformal homeomorphism in [9]. Our definition is
slightly stronger than theirs in this case (see [9] and [7]).
\enddemo

By the non-degeneracy of the contact structure of $M$, i.e., $d \eta \wedge
\eta \ne 0$ on $M$, there is an unique smooth vector field $T$ on $M$,
such that $T\lrcorner d\eta =0, \langle \eta, T \rangle =1$ on $M$.
$T$ is called the characteristic vector field for $\eta$.

Let $T^{1,0}M_0$ denote the subbundle $\{X-iJ_0X\,\vert\, X\in HM_0\}$ of
$\Bbb C \otimes TM_0$. Its elements are called $(1,0)$ vectors on $M_0$.
$T^{0,1}M_0\triangleq\overline {T^{1,0}}M_0$ is called $(0,1)$ tangent
bundle of $M_0$. Denote by $\wedge^{0,1}M_0$ the space of complex linear
functionals $\alpha$ on $\Bbb C\otimes HM$ so that $\alpha (Z)=0, \forall
Z\in T^{1,0}M_0$. An $\alpha \in \wedge^{0,1}M_0$ is called a $(0,1)$
form on $M_0$. Denote also $\overline {\wedge^{0,1}}M_0$ by $\wedge^{1,0}M_0$.

With two CR structures $M_0$ and $M_1$ on $M$ with the same orientation,
we associate a global section $\mu$ of $T^{1,0}M_0\otimes\wedge^{0,1}M_0$
as follows. Let $\overline W_0 \ne 0$ be a smooth $(0,1)$ vector field on
an open set $U \i M$ with respect to $M_0$, then $\mu$ is a section of
$T^{1,0}M_0\otimes\wedge^{0,1}M_0$ on $U$ so that $\overline W_1 =
\overline W_0 - \mu (\overline W_0 )$ is a (0,1) vector with respect to
$M_1$ on $U$. Let $\psi$ be a smooth $(1,0)$ form on $U$ with respect to
$M_0$ such that $\{\psi, \overline \psi\}$ is the dual basis to $\{W_0,
\overline W_0 \}$. With these conventions, $\mu = \nu W_0 \otimes \overline
\psi$ for a function $\nu$ on $U$. The tensor $\mu$ is globally well defined
and is called {\it the deformation tensor} of $M_1$ with respect to $M_0$. 
$|\mu |$ ($\triangleq |\nu |$ on $U$) is also a globally defined real valued
function. Since $M_0$ and $M_1$ have the same orientation, $|\mu |<1$
everywhere. 

\proclaim {Definition 2.2}If $f:M_1 \to M_0$ is a $C^1$ contact mapping which
preserves the orientation of $HM$, let $f^{-1}(M_0)$ be a new CR structure on
$M$ so that $T^{0,1}f^{-1}(M_0) =f^{-1}_*(T^{0,1}M_0)$. Define the Beltrami
tensor of $f$ by the deformation tensor of $f^{-1}(M_0)$ with respect to $M_1$.
\endproclaim

\demo {Remark} Locally, since
$$
f_*(\overline W_1) = \langle\psi , f_*(\overline W_1)\rangle W_0 +
\langle\overline\psi ,f_*(\overline W_1)\rangle\overline W_0 ,
\tag2.4
$$
we have
$$
\mu _f =\frac {\langle f_*\psi,\overline W_1\rangle }{\langle f_*\psi,
W_1\rangle } W_1\otimes \overline \psi_1 ,  \tag2.5
$$
where $\overline\psi_1\in\wedge^{0,1}M_1$ with $\langle\overline\psi_1,
\overline W_1\rangle = 1$. Since $f$ preserves the orientation of $HM$ and the
CR structures $M_0$ and $M_1$ have the same orientations, $\langle f_*\psi,
W_1\rangle \ne 0$ and $|\mu_f|<1$. Hence (2.5) and (2.6) below are
meaningful.
\enddemo

\proclaim {Theorem 2.3} If $f: M_1 \to M_0$ is a $C^1$ quasiconformal
homeomorphism and preserves the orientation of $HM$, then for $q \in M_1$, the
dilatation at the point $q$ is given by
$$
K(f)(q)=\frac {1+|\mu_f(q)|}{1-|\mu_f(q)|}. \tag2.6
$$

\noindent In particular, the maximal dilatation is 
$$
K(f)= \displaystyle\sup_{M_1} \frac {1+|\mu_f|}{1-|\mu_f|}=
\frac {1+\displaystyle\sup_{M_1}|\mu_f|}{1-\displaystyle\sup_{M_1}|\mu_f|}.
\tag2.7
$$
\endproclaim

The proof of this theorem is simple linear algebra and is the same
as the proof of an analogous fact on $\Bbb C$ (see [17]).

\vskip0.1truein
We now turn our attention to contact flows. First recall that
the non-degeneracy of the contact structure of $M$ shows that the mapping
$$
\imath : HM \to \text {Null}(T), \qquad X \mapsto X \lrcorner d\eta \tag2.8
$$ 
\noindent is a bundle isomorphism. Here the space 
$$
\text {Null}(T)=\{\omega \in \wedge ^1M\,\vert\,\langle \omega ,T
\rangle=0 \} \tag2.9
$$
is a real rank $2$ subbundle of $\wedge ^1M$. Denote the inverse of $\imath$
by $\sharp $.

Let $V$ be a vector field on a contact manifold $M$ which generates a
smooth flow of contact transformations. For such a vector field
$V$ the real valued function $u= \langle \eta , V \rangle$ is called
{\it the contact Hamiltonian function of $V$}.

\proclaim {Theorem 2.4 {\rm (i) }(Liebermann)} Suppose $M$ is a smooth compact
contact manifold with a smooth contact form $\eta$. If $V$ is a smooth
vector field which generates a flow of contact transformations of $M$, then 
$$
V=uT+\sharp ((Tu)\eta-du) ,\tag2.10
$$
here $u$ is the contact Hamiltonian of $V$. 

{\rm (ii)} Conversely, if $V$ is a vector field defined by (2.10) for a real
valued smooth function $u$ on $M$, then $V$ generates a flow of contact
transformations of $M$ and the Hamiltonian of $V$ is $u$.
\endproclaim

The part (i) is Th\'eor\`eme 3 in [15], a proof was given there. The  
sufficiency (ii) can be proved by straightforward computations.

\medskip 

On the $3$-sphere $S^3=\{ (w_1,w_2)\in \Bbb C^2 \,\vert \, |w_1|^2+
|w_2|^2=1\}$, the contact structure is defined by the contact form
$$
\eta =-\text {Im}(w_1d\overline w_1+w_2d \overline w_2) . \tag2.11
$$
The characteristic vector field for $\eta $ is
$$
T=-2\,\text {Im}(w_1{\partial \over \partial w_1}+w_2 {\partial \over \partial
w_2}) . \tag2.12
$$

Let $S_0^3$  be the sphere with the CR structure inherited
from the standard complex structure of $\Bbb C^2$. Let us denote
$$
\align
&W_0=\overline w_2 {\partial \over \partial w_1}-\overline w_1{\partial \over 
\partial w_2} , \tag2.13 \\
& \psi =w_2dw_1-w_1dw_2 . \tag2.14
\endalign
$$
Then $W_0, \overline W_0$ are $(1,0), (0,1)$ vector fields on $S_0^3$ 
respectively, and $\psi, \overline \psi$ are $(1,0), (0,1)$ forms on
$S_0^3$ respectively. Moreover $\{ W_0, \overline W_0, T\}$ is dual to $\{ \psi
, \overline \psi , \eta\}$. Direct computations yield the commutator relations
among these basis vectors of $\Bbb C \otimes TS^3$:
$$
[W_0,\overline W_0]=-iT, \qquad [T, W_0]=-2iW_0, \qquad [T,\overline W_0]
=2i\overline W_0 .
\tag2.15
$$ 

The vector fields $X\triangleq 2\text {Re}W_0,\, Y\triangleq -2\text {Im}W_0$
form a basis of the real contact space $HS^3$. We have
$$
[X,Y]=-2T, \qquad [X,T]=2Y, \qquad [Y, T]=-2X . \tag2.16
$$ 
The forms $\sigma \triangleq \text {Re}\psi ,\tau \triangleq \text {Im}\psi$
and $\eta$ form a basis of the cotangent space $\wedge^1S^3$.
The commutator relations (2.16) imply that $\imath (X)=2\tau, \imath
(Y)=-2\sigma$, or, equivalently, $\sharp (\tau )={1 \over 2}X, \sharp
(\sigma )=-{1 \over 2}Y$. So for any real valued function $u$ on $S^3$
$$
\sharp((Tu)\eta -du)=\sharp(-(Xu)\sigma-(Yu) \tau)=-{1\over 2}(Yu)X+{1\over 2}
(Xu)Y . 
$$
Hence we have proved the following corollary of Theorem 2.4.

\proclaim {Corollary 2.5} A vector field on $S^3$ generates a smooth
$1$-parameter group of contact transformations if and only if
$$
V=-{1 \over 2}(Yu)X+{1 \over 2}(Xu)Y+uT , \tag2.17
$$ 
or, equivalently,
$$
V=i(\overline W_0u)W_0-i(W_0u) \overline W_0+uT,
\tag2.18
$$ 
for a smooth real valued function $u$ on $S^3$. 
\endproclaim

\demo {Remark} An equivalent theorem in the setting of the 3-dimensional
Heisenberg group was given by Kor\'anyi and Reimann ([11], Theorem 5).
\enddemo

\heading
3. $S^1$-equivariant Quasiconformal Homeomorphisms
\endheading

Let $M$ be a smooth, compact 3-manifold. An $S^1$-action $\{ U_\phi\,\vert\,
\phi\in \Bbb R \mod 2\pi\}$ on $M$ is said to be free if no $U_\phi\ne
\text{id}$ has a fixed point. $M$ is called a regular contact manifold
if $M$ is contact and has a contact form $\eta$ so that the characteristic
vector field $T$ for $\eta$ generates a free $S^1$-action $\{ U_\phi\,\vert\,
\phi\in \Bbb R \mod 2\pi\}$ on $M$. Here $\phi$ is the parameter of the
contact flow. Obviously the action is transversal to the contact
structure. Let $\Sigma = M/S^1$ be the space of orbits. Then $\Sigma$ is
a smooth compact surface and the natural projection $p:M \to\Sigma$ is open
and smooth.

\proclaim {Theorem 3.1 (Boothby-Wang [3])} If $M$ is a regular contact
manifold, then

(i) $M$ is a principal fiber bundle over $\Sigma$ with structure
group $S^1$;

(ii) the contact structure $HM$ defines a connection in this bundle; and

(iii) $\Sigma$ has an oriented area form $\omega$ such that the structure
equation of the connection is given by $$d\eta = p^*\omega .$$
\endproclaim

Later we will simply call such a manifold $M$ {\it a contact circle bundle}.

A curve on a smooth compact manifold is said to be rectifiable if it is
rectifiable with respect to a (hence any) smooth Riemannian metric on the
manifold.

\proclaim{Lemma 3.2} Let $\gamma: I\to\Sigma$ be a rectifiable curve starting
at $q\in\Sigma$ with an interval $I= [0, l]\i\Bbb R$. $\tilde q\in p^{-1}(q)$.
Then there is a unique curve $\tilde\gamma: I\to M$ starting at $\tilde q$ so
that $p\circ\tilde\gamma =\gamma$, $\tilde\gamma$ is rectifiable, and the
tangent vectors at its regular points are in $HM$.
\endproclaim

The curve $\tilde\gamma$ is called the horizontal lift starting at $\tilde q$
of $\gamma$.

\demo{Proof} If $\gamma$ is $C^1$, the lemma follows from Proposition II 3.1
in [8]. The following is a modification of the proof given there.

By the local triviality of the circle bundle, we have a rectifiable curve
$\tilde\alpha :I\to M$ starting at $\tilde q$ so that $p\circ\tilde\alpha
=\gamma$. We construct an absolutely continuous function $\phi :I\to\Bbb R$
such that the curve given by
$$
\tilde\gamma (t) =U_{\phi(t)}(\tilde\alpha(t)), \qquad t\in I, \tag3.1
$$
satisfies the requirement. Note that if $T$ denotes the generator of the
circle action,
$$
\tilde\gamma^\prime(t)=\phi^\prime(t)\, T|_{{}_{\tilde\gamma(t)}}
+{U_{\phi(t)}}_*(\tilde\alpha^\prime(t)). \tag3.2
$$
This vector is in $HM$ if and only if
$$
0 = \langle\eta ,\tilde\gamma^\prime(t)\rangle =\phi^\prime(t) +\langle
\eta, {U_{\phi(t)}}_*(\tilde\alpha^\prime(t))\rangle. \tag3.3
$$
The expression on the right hand side of the ordinary differential equation
in the initial value problem
$$
\alignedat1
& {\phi}^\prime  = -\langle\eta, {U_{\phi}}_*(\tilde\alpha^\prime(t))
\rangle ,   \\
& \phi (0) = 0,
\endalignedat
\tag3.4
$$
is smooth in $\phi$ and $L^1$ in $t$. So, by Theorem II 3.5 in [18], (3.4) has
a unique solution $\phi$ on $I$ which is absolutely continuous. Then the curve
given by (3.1) with this $\phi$ is the horizontal lift starting at $\tilde q$
of $\gamma$.  \qed
\enddemo

Let $\Omega$ be a simply connected domain on $\Sigma$ with a rectifiable 
boundary $\gamma = \partial\Omega$. As an 1-chain $\gamma$ has an orientation
induced from that of $\Omega$ regarded as a 2-chain. For $q\in\gamma ,\tilde
q\in p^{-1}(q)$, let $\tilde\gamma$ be the horizontal lift of $\gamma$
starting at $\tilde q$. The end point of $\tilde\gamma$ is $U_\phi (\tilde q)$
for some $\phi\in [0, 2\pi )$. We call $\phi$ the {\it phase shift} from
$\tilde q$ to $U_\phi (\tilde q)$. The structure equation in Theorem 3.1
(iii) is the infinitesimal version of the following.

\proclaim {Proposition 3.3} The $\omega$-area of $\Omega$ satisfies
\hskip0.15truein $\dsize\int_\Omega\omega = -\phi\mod 2\pi$.
\endproclaim

\demo{Proof} Without loss of generality, we assume that $\Omega\subset
\subset\Omega^\prime$ for a simply connected open set $\Omega^\prime
\subset\Sigma$ where the bundle $M$ is trivial. That is,
$p^{-1}(\Omega^\prime )$ is $S^1$-equivariantly diffeomorphic to
$\Omega^\prime\times S^1$. Note $d\omega = 0$ on $\Sigma$, so $\omega
= d\alpha$ on $\Omega^\prime$ for some 1-form $\alpha$. Then
$$
\int_\Omega\omega =\int_\gamma\alpha =\int_{\tilde\gamma}p^*\alpha . \tag 3.5
$$
Here the first equality is due to the Stokes formula for rectifiable $\gamma$
which can be proved by exhausting $\Omega$ with $C^1$ bounded domains.
Notice the homology group $H_1(p^{-1}(\Omega^\prime))\cong \Bbb Z$. Let
$\beta$ be an $S^1$-fiber with the orientation given by $T$. Then regarded
as an 1-chain, $\beta$ generates $H_1(p^{-1}(\Omega^\prime))$.
If $\tilde\gamma_0$ is the oriented trajectory of $T$ from $\tilde q$
to $U_\phi(\tilde q)$, then $\tilde\gamma -\tilde\gamma_0$ is homologous
to $m\beta$ for some $m\in\Bbb Z$. Because
$$
\int_\beta p^*\alpha =\int_{p(\beta )}\alpha =0
$$
and
$$
d(\eta -p^*\alpha) = d\eta -p^*d\alpha =d\eta -p^*\omega =0, \tag3.6
$$ 
$$
\int_{\tilde\gamma -\tilde\gamma_0}\eta -p^*\alpha =\int_{m\beta}\eta-p^*
\alpha =\int_{m\beta}\eta =0 \mod 2\pi. \tag3.7
$$
Note also $\dsize\int_{\tilde\gamma}\eta =0$ since $\tilde\gamma$ is
Legendrian and $\dsize\int_{\tilde\gamma_0}p^*\alpha =\int_{p(\tilde\gamma_0)}
\alpha =0$. So (3.7) gives
$$
\int_{\tilde\gamma_0}\eta +\int_{\tilde\gamma}p^*\alpha = 0 \mod 2\pi,
$$
or, by (3.5),
$$
\int_\Omega\omega = -\int_{\tilde\gamma_0}\eta =-\phi \mod 2\pi. \qed
$$
\enddemo  

If we start with an oriented, rectifiable, Legendrian curve $\tilde
\gamma$ with the initial and end points on the same $S^1$-fiber, then
the closed curve $\gamma =p(\tilde\gamma)\i\Sigma$ may not bound a simply
connected domain, and $\tilde\gamma$ may not be a single-sheeted cover of
$\gamma$. However, when $\gamma$ represents the null element of $H_1(\Sigma )$
it is easy to see that Proposition 3.3 can be generalized to 

\proclaim {Corollary 3.4} If $p(\tilde\gamma )=\partial\Omega$ for some
2-chain $\Omega$ on $\Sigma$, the $\omega$-area of $\Omega$ has the same
value as the phase shift from the end point of $\tilde\gamma$ to its initial
point $(\text{mod } 2\pi)$.
\endproclaim 

A CR structure on $M$ is $S^1$-invariant if each $U_\phi$ in the
$S^1$-action is CR with respect to this CR structure. Assume $M_0$ is an
$S^1$-invariant CR manifold with the underlying regular contact manifold $M$,
then the CR structure induces a complex structure on the surface $\Sigma$ so
that $p: M\to\Sigma$ is CR. Equipped with this complex structure, $\Sigma$
becomes a Riemann surface $\Sigma_0$ and $T^{1,0}\Sigma_0 = p_*(T^{1,0}M_0)$.

Moreover, when the CR structure of $M_0$ is positively oriented with respect to
$\eta$, the area form $\omega$ and the complex structure on $\Sigma_0$
determine a Riemannian metric as follows. Let $J^\prime :T\Sigma_0\to
T\Sigma_0$ be the endomorphism which defines the complex structure on
$\Sigma_0$, then $\omega (X, J^\prime X) >0$ for nonzero $X\in T\Sigma_0$.
Then for $X, Y\in T\Sigma_0$, define a Riemannian metric by $\langle X,
Y\rangle =\omega (X,J^\prime Y)$. This Riemannian metric has the oriented
area form $\omega$ and induces the complex structure $J^\prime$ of $\Sigma_0$.
Still use $\Sigma_0$ to denote the corresponding Riemannian 2-manifold.

Conversely, if there is a Riemannian metric on $\Sigma$ whose oriented 
area form is $\omega$, we can lift the complex structure determined by this
Riemannian metric to an $S^1$-invariant CR structure on $M$ by declaring
$Z\in \Bbb C\otimes HM$ to be a $(1,0)$ tangent vector if $p_*(Z)\in T^{1,0}
\Sigma$. This CR structure is positively oriented with respect to $\eta$.

A homeomorphism $f:M\to M$ is said $S^1$-equivariant if the diagram
$$
\CD 
 M @> f >>M \\
@V U_\phi VV @VV U_\phi V \\ 
 M @> f >>M
\endCD
\tag 3.8
$$
commutes for each $\phi$. Such a homeomorphism will induce a quotient
homeomorphism $F:\Sigma\to\Sigma$ so that the diagram
$$
\CD
 M @> f >> M \\
@V p VV @VV p V \\
 \Sigma @> F >> \Sigma 
\endCD
\tag 3.9
$$
commutes.

Assume $M_1$ is another $S^1$-invariant CR manifold with the underlying
contact manifold $M$. The corresponding quotient surface is $\Sigma_1
=M_1/S^1$ which has the area form $\omega$ too and the complex
structure induced from the CR structure on $M_1$.

\proclaim {Theorem 3.5} Let $M \overset{p}\to{\to}\Sigma$ be a contact
circle bundle. Assume $M_1, M_0$ are two $S^1$-invariant CR manifolds
with the same underlying contact manifold $M$, and $f: M_1\to M_0$ is an
$S^1$-equivariant quasiconformal homeomorphism. Then the quotient map
$F: \Sigma_1\to\Sigma_0$ is a quasiconformal homeomorphism in the classical
sense and $F$ preserves $\omega$-area. Moreover $K(F)=K(f)$.
\endproclaim

\demo {Proof} Choose a region $R$ on the Riemann surface $\Sigma_1$
corresponding to a rectangle in a conformal coordinate system. Let $\Gamma =
\{\gamma \}$ be the family of all longest straight line segments in $R$ which
are parallel to a fixed side of $R$. Lifting each $\gamma\i\Gamma$ to $M_1$
horizontally, we obtain a contact fibration $p^{-1}(\Gamma )=\{ $all
Legendrian lifts of $\gamma\,\vert\,\gamma\in \Gamma\}$ of $p^{-1}(R)$. Let
$\Gamma _1\i\Gamma$
consist of lines $\gamma$ so that $f$ is absolutely continuous along
a lift of $\gamma$. $S^1$-equivariance tells us
if $\gamma \in \Gamma _1$, then along each lift of $\gamma $, $f$ is
absolutely continuous. Therefore if $\gamma \in \Gamma _1$, then $F$ is
absolutely continuous along it. By the ACL property of $f$, $p^{-1}(\Gamma
\setminus\Gamma _1)$ is of measure zero. Therefore, $F$ is absolutely
continuous along almost every straight line segment $\gamma \in \Gamma$.
Since $R$ is arbitrary, $F$ is ACL.

If $f$ is differentiable at a point $\tilde q$, $F$ is differentiable at
$ q=p(\tilde q)$. Hence $F$ is differentiable almost everywhere on $\Sigma$
since so is $f$ on $M$. The bounded distortion inequality for $f$ at $\tilde q$ 
implies that for $F$ with the same dilatation at $q$ since $p$ is CR. So
$F$ is a quasiconformal homeomorphism of $\Sigma$ and $K(f)=K(F)$.

For $q\in \Sigma_1$, let $D_r$ be a disc with radius $r$ centered at $q$,
for each positive small $r$. ACL regularity and $S^1$-equivariance of
$f$ implies that $F$ is absolutely continuous along almost all circles
$\partial D_r$, and $f$ is absolutely continuous along all lifts of these
circles. For those discs $D_r$ along whose boundary $F$ is absolutely
continuous (equivalently, $f$ is absolutely continuous along each lift of
$\partial D_r$), $F(\partial D_r)$ is rectifiable. Hence Proposition 3.3 
is valid for both such $D_r$ and the corresponding $F(D_r)$. Then
$S^1$-equivariance of $f$ and Proposition 3.3 show that $F$ preserves
the $\omega$-area of almost all discs $D_r$, hence of all discs.
So $F$ preserves the $\omega$-area for $q$ is arbitrary. \qed
\enddemo

When $\Sigma$ is simply connected and $F: \Sigma_1\to\Sigma_0$ is $C^1$,
we have the following converse to Theorem 3.5.

\proclaim{Theorem 3.6} Let $M\overset{p}\to{\to}\Sigma$ be a compact
contact circle bundle with $\Sigma$ homeomorphic to $S^2$. For $j=0,1$,
let $\Sigma_j$ be a Riemannian 2-manifold obtained by assigning to $\Sigma$ a
Riemannian metric whose area form is $\omega$; let $M_j$ be an $S^1$-invariant
CR manifold obtained by endowing $M$ with the CR structure such that
$p: M_j\to\Sigma_j$ is CR. Assume $F:\Sigma_1\to\Sigma_0$ is a $C^1$
quasiconformal homeomorphism which preserves $\omega$-area. Then there
exists an equivariant quasiconformal homeomorphism $f:M_1\to M_0$
such that $p\circ f=F\circ p$ and $K(F)=K(f)$.
\endproclaim

\demo{Proof} Fix a point $q_0 \in \Sigma_1$ and a points ${\tilde q}_0 \in
p^{-1}(q_0)$. Define $f(\tilde q_0)$ to be any point in the fiber $p^{-1}
(F(q_0))$. For any other $\tilde q \in M_1$, connect $\tilde q_0$ and
$\tilde q$ by a $C^1$ Legendrian curve $\tilde \gamma$. We can always
do that by a theorem of Chow [5]. Project $\tilde \gamma$ onto a curve
$\gamma \i \Sigma_1$, then map it by $F$ onto the $C^1$ curve
$F(\gamma )\i \Sigma_0$. We define $f(\tilde q)$ by the end point of
the unique horizontal lift of $F(\gamma )$ starting at $f(\tilde q_0)$.  

Assume $\tilde \gamma _1$ is another $C^1$ Legendrian curve connecting
$\tilde q_0$ and $\tilde q$, and $\gamma _1$ is its projection. Since $\Sigma$
is simply connected, the 1-chain $\gamma_1 -\gamma =\partial\Omega$ for
some 2-chain $\Omega\i\Sigma_1$. Corollary 3.4 says that the $\omega$-area
of $\Omega$ is zero$\mod 2\pi$, whence the same holds for the $\omega$-area of
$F(\Omega )$ since $F$ preserves $\omega$-area. By Proposition 3.3, the
horizontal lift of $F(\gamma )$ and $F(\gamma _1)$ initiated at $f(\tilde q_0)$
have the same end points. Therefore the mapping $f$ is well-defined. 

The map $f$ defined above is a $C^1$ contact homeomorphism, by the $C^1$
dependence of the horizontal lift of $F(\gamma)$ on $F(\gamma)$ which follows
the theorem in Appendix 1 of [8]. $f$ is also $S^1$-equivariant by an
argument similar to the one given in the last paragraph based on Corollary
3.4. Its bounded distortion inequality follows from that of $F$, and $f$, $F$ 
share the same value of dilatation since the $S^1$-action is CR. \qed
\enddemo

\demo {Remark} (1) The lift $f$ of $F$ constructed in the proof is unique up to
composition with $U_\phi$ for some $\phi$.

(2) When the base space $\Sigma$ is not simply connected, a $C^1$
homeomorphism $F$ on $\Sigma$ preserving $\omega$ can be lifted to a
differentiable homeomorphism $f$ whose differential preserves the contact
structure if and only if the monodromy representation of $\pi_1 (\Sigma)$
in $S^1$ induced by $F$ is trivial. In this case, the construction
of $f$ in the above proof applies. When $\Sigma$ is homeomorphic to $S^2$, this
obstruction to lifting does not exist. 
\enddemo

When $M=S^3=\{|w_1|^2+|w_2|^2=1\}\i {\Bbb C}^2$ and the circle action
is given by
$$
U_\phi : (w_1, w_2)\to (e^{i\phi}w_1, e^{i\phi}w_2), \tag3.10
$$
we have the Hopf fibration $S^1\to S^3 \to S^2$ of the 3-sphere. The projection
is given by
$$
p:S^3\to S^2, \hskip0.1truein (w_1, w_2) \mapsto {w_2 \over w_1}. \tag 3.11
$$
On $S^2$ the standard spherical metric is
$$
ds=\frac {2\,|dz|}{1+|z|^2} \tag 3.12
$$
and $\dsize\omega_0={4\,dx\wedge dy\over (1+|z|^2)^2}$ is the spherical area form,
where $z=x+yi$. Let $\eta$ be the contact form of $S^3$ given by (2.11). Then
direct computations prove

\proclaim{Proposition 3.7} We have $\dsize \qquad d\eta =p^*({1\over 2}\,
\omega_0).$
\endproclaim

\heading
4. Equivariantly Extremal Quasiconformal Homeomorphisms on $S^3$
\endheading

Here an equivariantly extremal quasiconformal homeomorphism refers
to an equivariant quasiconformal homeomorphism with the least maximal
dilatation among all equivariant homeomorphisms.

Given two smooth Riemannian metrics on $S^2$ which share the spherical
area form, we lift the complex structures they determine to two smooth
$S^1$-invariant CR structures on $S^3$ so that the projection $p$ in (3.11)
is CR. By results in the last section, if an extremal area-preserving
quasiconformal homeomorphism on $S^2$ between these two Riemannian structures
is $C^1$, then an $S^1$-equivariant lift of this homeomorphism is an
$S^1$-equivariant extremal quasiconformal homeomorphism on $S^3$ between two
lifted CR structures. This is the guideline for this section.   

The spherical metric on the unit Euclidean sphere $S_0^2$ is given by (3.12),
or, equivalently,
$$
ds_0^2=d\theta ^2+\sin^2\theta d\phi ^2 , \tag4.1
$$
where $(\theta ,\phi)$ are the spherical coordinates $(0\le \theta \le \pi ,
0\le\phi <2\pi)$. Let $\lambda$ be a real valued smooth function on $S^2$
satisfying $1\le\lambda\le \Lambda $ on $S^2$, $\lambda =1$ near the poles
where $\theta = 0, \pi$, $\lambda$ attains its maximal value $\Lambda >1$ at
each point of the equator $E=\{ \theta ={\pi \over 2} \}$, and $\lambda
< \Lambda$ elsewhere. Define a new metric on
$S^2$ by 
$$
ds_1^2=\lambda^2d\theta ^2+\frac {\sin^2\theta }{\lambda^2}{d\phi^2}. \tag4.2
$$
$S^2$ equipped with the metric (4.2) is denoted by $S_1^2$. The metric on
$S_1^2$ is obtained from the metric on $S_0^2$ by stretching in the meridian 
direction by the factor $\lambda$ and shrinking in the parallel direction
by the same factor. $\text{id}_{S^2}: S_1^2\to S_0^2$ is quasiconformal
with maximal dilatation $\Lambda^2$ which occurs along the equator. Obviously,
$S_0^2$ and $S_1^2$ have the area element $\sin\theta d\theta d\phi$.

A Jordan curve  divides the sphere into two components. If these components
have equal area, we call the curve area-halving curve. An area-halving curve
on $S_0^2$ is also an area-halving curve on $S_1^2$. Let us give a folk lemma
first. It is a very special case of isoperimetric property on surfaces
(Burago and Zalgaller [4], Theorem 2.2.1.). Our proof is very simple and
intuitive.

\proclaim {Lemma 4.1} The great circles on $S_0^2$ are the shortest
area-halving curves.
\endproclaim

\demo {Proof} Any two area-halving curves on $S_0^2$ must intersect each
other. Hence an area-halving curve intersects its antipodal image,
and we conclude that an area-halving curve contains a pair of antipodal
points. But the semi-great circles are the geodesics to connect two antipodal
points. Therefore a Jordan curve is a shortest area-halving curve if and
only if it is a great circle. \qed
\enddemo

\smallskip Therefore the length of a shortest area-halving curve on $S_0^2$
is $2\pi$. The construction of $ds_1^2$ shows that on $S^2_1$ the equator
is the unique shortest area-halving curve and its length is ${2\pi / \Lambda }$.

\proclaim {Proposition 4.2} The identity map $\text{id}_{S^2}: S_1^2 \to S_0^2$
has the least maximal dilatation among all area-preserving quasiconformal
homeomorphism from $S^2_1$ to $S^2_0$.
\endproclaim

\demo {Proof} Divide the equator $E=\{\theta ={\pi \over 2} \}\i
S_1^2$ by ordered points $q_1,q_2,...,q_n$ ($q_{n+1}=q_1$) into small
subarcs. Let the $\phi$-coordinate of $q_j$ be $\phi_j$. For $1 \le j \le n$
and small $\delta > 0$, form a quadrilateral $Q$ given by
${\pi\over 2}-\delta\le\theta\le {\pi \over 2}$, $\phi _j\le\phi\le
\phi_{j+1}]$. Then the four vertices of $Q$ are $q_j$, $q_{j+1}$,
$p_{j+1}$ and ${p_j}$ for some points $p_{j+1}$ and $p_j$ on a parallel.
Recall the module of the quadrilateral $Q$ is defined by
$$
\text{Mod}(Q)=\displaystyle \sup_{\varrho \in A(Q)} \frac {(\displaystyle
\inf_{\gamma\in\Gamma _Q}\int_\gamma\varrho )^2}{\int_Q\varrho^2}
=\displaystyle\inf_{\varrho\in A(Q)}\frac {\int_Q\varrho^2}{(\displaystyle
\inf_{\gamma\in\Gamma_Q^{\prime}}\int_\gamma\varrho )^2}, \tag4.3
$$
where $A(Q)=\{\varrho \ge 0 \, \vert$ $\varrho $ is Borel-measurable on $Q$,
$0< \int_Q \varrho^2<+\infty \}$ is the set of allowable measures, $\Gamma_Q$
is the family of rectifiable curves in $Q$ connecting the sides $q_jq_{j+1}$,
$p_jp_{j+1}$, and $\Gamma_Q^{\prime}$ is the family of rectifiable curves
in $Q$ connecting the sides $p_jq_j$, $p_{j+1}q_{j+1}$. In particular
$$
\frac {(\displaystyle\inf_{\gamma\in\Gamma_Q} \int_\gamma 1 )^2}{\text {Area}
(Q)}\le\text{Mod}(Q) \le \frac {\text {Area}(Q)}{(\displaystyle \inf_{\gamma
\in \Gamma_Q^{\prime}} \int_\gamma 1)^2}. \tag4.4
$$ 
Similar definitions and inequalities hold for the quadrilateral $F(Q)$.
If a homeomorphism $F: S^2_1\to S^2_0$ is $K$-quasiconformal,
$$
\text{Mod}(Q) \le K\text{Mod}(F(Q)).
$$
Combining this with (4.4) for both $Q$ and $F(Q)$, we have 
$$
\frac {(\displaystyle\inf_{\gamma\in\Gamma_Q} \int_\gamma 1)^2}{\text {Area}
(Q)}\le K \frac {\text {Area}(F(Q))}{(\displaystyle \inf_{\gamma \in 
\Gamma_{F(Q)}^{\prime}} \text {length}(\gamma ) )^2}. \tag4.5
$$ 
Denote $d =d(\delta ) \triangleq \displaystyle \inf_{\gamma \in \Gamma_{F(Q)}
^{\prime}} \text {length}(\gamma )$. This is the distance between the side
$F(p_j)F(q_j)$ and the opposite side $F(p_{j+1})F(q_{j+1})$ of $F(Q)$. Hence 
$$
\displaystyle \lim_{\delta \to 0} d(\delta )=d_0(F(q_j), F(q_{j+1})),
$$
here $d_0$ is the distance on $S_0^2$. Since $F$ preserves the area,
$$
\text {Area}(F(Q))=\text {Area}(Q)=\int_{{\pi \over 2}-\delta }^{\pi \over 2}
\int_{\phi_j}^{\phi_{j+1}} \sin\theta d\theta d\phi =|\phi_{j+1}-\phi_j|\sin
\delta .
$$
Then (4.5) becomes
$$
d \displaystyle\inf_{\gamma\in\Gamma_Q} \int_\gamma 1 \le \sqrt K |\phi
_{j+1}-\phi_j|\sin\delta .
$$
For any $\epsilon > 0$, there exists a $\gamma \in \Gamma_Q$ such that
$$
d\, ({1 \over \delta}\int_\gamma 1) \le \sqrt K |\phi_{j+1}-\phi_j| \frac
{\sin \delta}{\delta } +\epsilon .
$$
Letting $\delta \to 0$,
$$
d_0(F(q_{j+1}), F(q_j))\Lambda \le \sqrt K |\phi_{j+1}-\phi_j|+\epsilon .
$$
Letting $\epsilon \to 0$, and summing over all $j$
$$
\Lambda \displaystyle \sum_{j=0}^nd_0(F(q_{j+1}),F(q_j)) \le \sqrt K
\displaystyle \sum_{j=0}^n |\phi_{j+1}-\phi_j| = 2\pi\sqrt K .
$$
By the arbitrariness of the partition of $E$, hence of the corresponding
partition of $F(E)$, we conclude that the area-halving curve $F(E)$ on $S_0^3$
is rectifiable and
$$
2\pi\sqrt K \ge\Lambda\text{ length}(F(E))\ge 2\pi\Lambda ,
$$
by Lemma 4.1. Therefore, $K(F)=K \ge \Lambda^2= K(\text{id}_{S^2})$. \qed
\enddemo

The Riemannian metric on $S^2_1$ given by (4.2) can be written as  
$$
ds^2_1=\frac {(\lambda^2 +1)^2}{\lambda^2 (|z|^2+1)^2}\,\big| dz+\frac
{\lambda^2 -1}{\lambda^2 +1}{z\over {\overline z}}d\overline z \big|^2.\tag4.7
$$
Then on $S^2_1$, the (0,1) tangent space is spanned by
$$
\frac {\partial }{\partial \overline z}- \frac {\lambda^2 -1}{\lambda^2 +1}
{z \over {\overline z}} \frac {\partial }{\partial z}, \tag4.8
$$
which is annihilated by the (1,0) form 
$$
dz+\frac {\lambda^2 -1}{\lambda^2 +1}{z \over {\overline z}}d\overline z.
$$

Denote $\tilde\lambda = \lambda\circ p$ and $\overline W_1=\overline W_0 -\nu
W_0$, where $W_0$ is given by (2.13) and $$\nu =\frac {{\tilde\lambda}^2 -1}
{{\tilde\lambda}^2 +1} \frac {w_1w_2} {\overline w_1 \overline w_2}.
\tag4.9$$ 
Then direct computations give
$$
p_*(-{\overline w_1}^2 {\overline W}_1)=\frac {\partial }{\partial \overline
z}- \frac {\lambda^2 -1}{\lambda^2 +1} {z \over {\overline z}}
\frac {\partial }{\partial z}. \tag4.10
$$ 

Use $S_1^3$ to denote $S^3$ equipped with the CR structure whose
(0,1) vector space is spanned by $\overline W_1$. By Theorem 3.5, 3.6, 3.7 and
Proposition 4.2, we have proved

\proclaim {Theorem 4.3} With above notation, $\text{id}_{S^3}: S_1^3 \to S_0^3$
is an equivariantly extremal quasiconformal homeomorphism, namely, it
has the least maximal dilatation among all equivariant quasiconformal
homeomorphism from $S_1^3$ to $S^3_0$.
\endproclaim

\demo {Remarks}
(1) The dilatation of $\text{id}_{S^3}:S_1^3 \to S_0^3$ attains its
maximum on the covering of the equator $E \i S^2$, i.e., the Clifford torus
$$
T_C=\{(w_1,w_2)\, |\, |w_1|^2=|w_2|^2={1 \over 2}\}
$$
and its maximal value is $\Lambda^2$.

(2) $\text{id}_{S^3}:S_1^3 \to S_0^3$ is not the only extremal extremal
$S^1$-equivariant quasiconformal homeomorphism. Any small $S^1$-equivariant
perturbation of $\text{id}_{S^3}$ away from $T_C$ will give another extremal
mapping.
\enddemo

\heading
5. Variation of the Conformal Distortion
\endheading

As before, we denote the 3-sphere endowed with the canonical CR structure
by $S_0^3$. Assume $S_1^3$ is the 3-sphere endowed with a new smooth,
strongly pseudoconvex CR structure whose $(0,1)$ tangent space is spanned by
$\overline W_1=\overline W_0-\mu (\overline W_0)$, where $\mu =\nu
W_0\otimes\overline\psi$ is a global section of $T^{1,0}S_0^3\otimes\wedge^{0,1}
S_0^3$ for a smooth function $\nu$ with $|\nu |<1$ on $S^3$.

Let $g_s$ be a flow of contact transformations generated by a vector
field $V$ with Hamiltonian function $u$. Then the maximal dilatation 
of $g_s:S_1^3 \to S_0^3$, by Theorem 2.3, is measured by the magnitude
of the Beltrami tensor $\mu_{g_s}$.

In this section we will give an asymptotic formula for $|\mu_{g_s}|$
as $s\to 0$ up to the first order for a general CR structure on $S_1^3$
and then up to the second order when the CR structure on $S_1^3$ is
$S^1$-invariant and the first variation vanishes.

According to (2.5) 
$$
|\mu_{g_s}|=\Big|\frac{\langle g_s^*\psi, \overline W_1\rangle}{\langle g_s^*
\psi, W_1\rangle}\Big| =\Big|\frac{\nu_s -\nu}{1-\overline\nu\nu_s}\Big|,
\tag5.1
$$
where 
$$
\alignedat1
\nu_s & \triangleq{\langle g_s^*\psi ,\overline W_0\rangle\over\langle
g_s^*\psi ,W_0 \rangle}  \\
& = \frac {\langle L_V\psi ,\overline W_0 \rangle s+{1 \over 2}\langle L_VL_V
\psi ,\overline W_0 \rangle s^2+\Cal{O} (s^3)}{1+\langle L_V\psi ,W_0\rangle s+
\Cal{O} (s^2)}  \\
& =\langle L_V\psi ,\overline W_0 \rangle s+\left( {1 \over 2}
\langle L_VL_V\psi , \overline W_0 \rangle -\langle L_V\psi , \overline W_0
\rangle \langle L_V\psi , W_0\rangle \right) s^2+\Cal{O} (s^3) \\
& \triangleq as+bs^2+\Cal{O} (s^3),
\endalignedat
\tag5.2
$$
for small $s\in \Bbb R$. Then on the set where $\nu \ne 0$,
$$
\alignedat1
|\mu_{g_s}| = & \, |(\nu -\nu_s)\left(1+\overline\nu\nu_s+{\overline\nu }^2
{\nu_s}^2 +\Cal{O} (s^3) \right) |  \\ 
= & \, |\nu |-\frac{1-|\nu |^2}{|\nu |}\text{Re}(\overline\nu a)s+\frac{1-
|\nu |^2}{2|\nu |}\bigl((1-|\nu |^2)|a|^2-2\text{Re}({\overline\nu}^2a^2)
-2\text{Re} (\overline \nu b)\bigr)s^2 \\ 
&  +\Cal{O} (s^3) .
\endalignedat
\tag5.3
$$
 
Now we compute the coefficients appearing in (5.2) and (5.3).
$$
\alignedat1
L_V\overline W_0 & =[V, \overline W_0]  \\
& =[i(\overline W_0u)W_0-i(W_0u) \overline W_0+uT, \overline W_0] \qquad\text
{by (2.18)}  \\
& =-i({\overline W_0}^2u)W_0+i(\overline W_0W_0u+2u)\overline W_0 ,
\endalignedat
\tag5.4
$$
and so
$$L_VW_0 =i(W_0^2u) \overline W_0-i(W_0 \overline W_0u+2u)W_0 . \tag5.5 $$
Hence
$$
a =\langle L_V\psi ,\overline W_0\rangle =V\langle \psi ,\overline W_0
\rangle -\langle \psi ,L_V \overline W_0\rangle =i({\overline W_0}^2u) .
\tag5.6
$$

Combining (5.3) with (5.6), we have proved the following proposition about the
first variation of the absolute value of Beltrami tensor.

\proclaim {Proposition 5.1} If $g_s:S_1^3 \to S_0^3$ is a flow of contact
transformations generated by a vector field with Hamiltonian $u$, then
for small $s \in \Bbb R$
$$
\alignat2
|\mu_{g_s}| & =|\nu |+\frac {1-|\nu |^2}{|\nu |}\text {Im}(\overline \nu
{\overline W_0}^2u)s+\Cal{O} (s^2) & \qquad \text{ where $\nu \ne 0$; and}
\tag5.7 \\
|\mu_{g_s}| & =|{\overline W_0}^2u| \cdot |s|+\Cal{O} (s^2) &
\qquad\text{ where $\nu = 0$.} \tag5.8
\endalignat
$$
\endproclaim

We will go on to compute the second order term in (5.2) and (5.3). By (5.5)
$$
\alignedat1
\langle L_V\psi ,W_0\rangle & =V\langle\psi ,W_0\rangle -\langle\psi ,L_V W_0
\rangle  \\
& =i(W_0\overline W_0u+2u),
\endalignedat
\tag5.9
$$
$$
\alignat2
\langle L_VL_V\psi ,\overline W_0\rangle= & V\langle L_V\psi ,
\overline W_0 \rangle -\langle L_V\psi ,L_V\overline W_0 \rangle  & \\
= & \Bigl( i(\overline W_0u)W_0-i(W_0u)\overline W_0+uT \Bigr)
(i{\overline W_0}^2u) & \text {by (2.18)} \\
& -\langle L_V\psi ,-i({\overline W_0}^2u)W_0+i(\overline W_0W_0u+2u)
\overline W_0 \rangle  & \text {by (5.4)} , \\
= & -(\overline W_0u)(W_0{\overline W_0}^2u)+(W_0u)({\overline W_0}^3u)
& \tag5.10 \\
& +iu(T{\overline W_0}^2u)-({\overline W_0}^2u)([W_0, \overline W_0]u) & \text
{by (5.6),(5.9)},\\
= & -(\overline W_0u)(W_0{\overline W_0}^2u)+(W_0u)({\overline W_0}^3u) & \\
&+iu(T{\overline W_0}^2u) +i({\overline W_0}^2u)(Tu) , & \text {by (2.15)}.
\endalignat
$$

So we finally get the expression of $b$ in (5.2).
$$
\alignedat1
b= & {1 \over 2}\langle L_VL_V\psi ,\overline W_0 \rangle -\langle L_V\psi ,
\overline W_0 \rangle \langle L_V\psi ,W_0\rangle  \\
= & -{1 \over 2}(\overline W_0u)(W_0{\overline W_0}^2u)+{1 \over 2}(W_0u)
({\overline W_0}^3u)+{1 \over 2}iu(T{\overline W_0}^ 2u) \\
& +{1 \over 2}i(\overline W_0^2u)(Tu)+(\overline W_0^ 2u)(W_0\overline
W_0u)+ 2(\overline W_0^2u)u.
\endalignedat
\tag5.11
$$

If on the set where $\mu \ne 0$, $\text{Im}(\overline \nu{\overline W_0}^2u)=0$,
i.e., the first variation of the absolute value of Beltrami tensor vanishes,
then Proposition 5.1 is not enough to analyse the behavior of the pertubation.
We will need to study the second variation of $|\mu_{g_s}|$ in this case.

Next we will compute the second order term in (5.3) on the set where 
$$
\text {Im}(\overline \nu{\overline W_0}^2u)=0 \qquad\text{and}\qquad\nu\ne 0
\tag5.12
$$
holds. Note one term in the second order coefficient in (5.3) is 
$$
\alignedat1
2\text {Re}(\overline \nu b) = & \text { Re}\left( -\overline \nu(\overline
W_0u) (W_0{\overline W_0}^2u)+\overline \nu(W_0u)({\overline W_0}^3u)
+2\overline \nu ({\overline W_0}^ 2u)(W_0\overline W_0u)\right) \\
& +\text {Re} \left( iu\overline \nu(T{\overline W_0}^2u)+4\overline \nu
({\overline W_0}^2u)\right) \\
& +\text{Re}\left( (i\overline \nu {\overline W_0}^2u) (Tu)\right) \\
\triangleq & \, I_1+I_2+I_3.
\endalignedat
\tag5.13
$$

To simplifiy $I_1$, let $\dsize c=\frac{\overline \nu ({\overline W_0}^2u)}{|\nu
|^2}$. With the assumption (5.12), $c$ is real valued.
$$
\alignedat1
I_1 & =\text{Re}\left(-\overline \nu (\overline W_0u)W_0(\nu c)+\overline
\nu (W_0u) \overline W_0 (\nu c)+2\overline \nu ({\overline W_0}^2u)
(W_0\overline W_0u)\right) \\
& =\overline \nu ({\overline W_0}^2u)\left( \Delta u+ \text{Re}\Bigl( {1 \over
\nu} (\overline W_0\nu)(W_0u)-{1 \over \nu}(W_0\nu )(\overline W_0u)
\Bigr)\right) , 
\endalignedat
\tag5.14
$$
where $\Delta u=(W_0\overline W_0+\overline W_0W_0)u$. 

For simplicity and for later applications, we will assume in the rest of this
section that the CR structure of $S_1^3$ is $S^1$-invariant. Then
$S^1$-invariance of the CR structure on $S_1^3$ implies that $L_T(\overline
W_0-\nu W_0)$ is a multiple of $\overline W_0-\nu W_0$. But
$$
\alignedat1
L_T(\overline W_0-\nu W_0) & =[T, \overline W_0-\nu W_0]  \\
& = 2i\overline W_0+(2i\nu -T\nu)W_0, \hskip1.5truein\text{by (2.15)}.
\endalignedat
\tag5.15
$$
Therefore, we have proved
\proclaim{Proposition 5.2} On $S^3$, $\mu =\nu W_0\otimes\overline\psi$ defines
an invariant CR structure if and only if
$$
L_T\mu = 4i\mu \qquad\text{or}\qquad T\nu =4i\nu .\tag5.16
$$
\endproclaim

With this simple fact, we have
$$
\alignedat2
I_2 & =\text{Re}\left( iuT(\overline\nu {\overline W_0}^2u)-iu(T\overline\nu)
{\overline W_0}^2u +4\overline\nu({\overline W_0}^2u)u \right)  &\qquad & \\
& =uT\left( \text{Re}(i\overline\nu{\overline W_0}^2u)\right) +\text{Re} \left(
-4u\overline\nu {\overline W_0}^2u+4u\overline\nu{\overline W_0}^2u \right),
& \qquad & \text{by (5.15)} , \\
& =0, & \qquad & \text{by (5.12)}.
\endalignedat
\tag5.17
$$

Obviously $I_3=0$ by (5.12). Combining this with (5.3), (5.6), (5.13) and
(5.14), we obtain

\proclaim {Proposition 5.3} If the smooth CR structure on $S_1^3$ is
$S^1$-invariant, the Beltrami tensor of $g_s:S_1^3 \to S_0^3$ satisfies
$$
\alignedat1
|\mu_{g_s}|= & |\nu |+{1-|\nu |^2 \over 2|\nu |} \Big\lbrace (1+|\nu |^2)| 
{\overline W_0}^2 u|^2-(\overline\nu{\overline W_0}^2u) \Big\lbrack \Delta u \\ 
& +\text{Re}\Bigl({1 \over \nu}(\overline W_0\nu)(W_0u)-{1 \over \nu}(W_0\nu )
(\overline W_0u)\Bigr) \Big\rbrack \Big\rbrace s^2+\Cal{O} (s^3),
\endalignedat
\tag5.18
$$
for small $s \in \Bbb R$ on the set where $\nu \ne 0$ and $\text{Im}
(\overline\nu{\overline W_0}^2u)=0$.
\endproclaim 

\heading
6. Symmetry Breaking
\endheading

In this section, we will use a contact perturbation of the equivariantly
extremal quasiconformal homeomorphism $\text{id}_{S^3}:S_1^3 \to
S_0^3$ constructed in Section 4 to show $\text{id}_{S^3}$ is not
extremal among all quasiconformal homeomorphisms between $S_1^3$ and $S_0^3$.
Namely, we will construct a nonequivariant quasiconformal homeomorphism
near $\text{id}_{S^3}$ with smaller maximal dilatation. That will prove
the following

\proclaim {Theorem 6.1} With $S_1^3$, $S_0^3$ denoting the $S^1$-invariant
CR manifolds constructed in section 4, no extremal quasiconformal
homeomorphism between $S_1^3$ and $S_0^3$ is equivariant. 
\endproclaim

\smallskip We call this phenomenon a symmetry breaking of the extremal
quasiconformal homeomorphism between CR structures on $S^3$.

\demo {Proof} Assume an extremal quasiconformal homeomorphism $f: S_1^3
\to S_0^3$ is equivariant. By Theorem 4.3, $K(f)=K(\text{id})$. We shall
construct a contact flow $g_s$ with a Hamiltonian $u$ which satisfies 
$$
\alignat2
& \text{Im}(\overline \nu {\overline W_0}^2u)=0, \qquad & \text{on } S^3,
\tag6.1 \\
& (1+|\nu |^2)|{\overline W_0}^2u|^2-(\overline \nu {\overline W_0}^2u)\Delta u
<0, \qquad & \text{on the torus } T_C. \tag6.2
\endalignat
$$
Here (6.1), by Proposition 5.1, makes the first variation of the absolute
value of Beltrami tensor of $g_s: S_1^3 \to S_0^3$ zero, and Proposition
5.3 applies. Direct computations show that $W_0\nu =\overline W_0 \nu=0$ on
$T_C$. So (6.2) gives that the the second order term in (5.18) is negative.
This will contradict the extremality of $f$, since $K(g_s) < K(f)$ for small
$s\in\Bbb R$.

For (6.2), we consider the equation
$$
(1+|\nu |^2)W_0^2u-\overline \nu\Delta u=-W_0^2u 
$$
on $T_C$. By (4.9) this is equivalent to
$$
\Delta u-H{w_1w_2 \over \overline w_1 \overline w_2} W_0^2u=0 \tag6.3
$$
on $T_C$, here $H$ is the constant value of ${\dsize 2+|\nu |^2\over \dsize |\nu |}$ on $T_C$.
Hence to satisfy (6.1), (6.2), it suffices to find $u$ satisfying the system
$$
\cases
\dsize\Delta u-H\,\text{Re}\left( {w_1w_2 \over \overline w_1 \overline w_2}
W_0^2u \right) =0, \qquad\qquad & \text{on } T_C,\\
\dsize\text{Re}\left( {w_1w_2 \over \overline w_1 \overline w_2}W_0^2u\right)
\ne 0, \qquad\qquad & \text{on } T_C,\\
\dsize\text{Im}\left( {w_1w_2\over\overline w_1\overline w_2}W_0^2u\right)=0,
\qquad\qquad & \text{on } S^3.
\endcases
\tag6.4
$$
If $u$ is independent of $w_2$, the system (6.4) is simplified to
$$
\cases
\dsize {\partial^2u \over\partial w_1\partial \overline w_1}- \text{Re}\left(
2w_1{\partial u \over \partial w_1}+Hw_1^2 {\partial^2u \over \partial w_1^2}
\right)=0, \qquad & \text{when } |w_1|^2={1\over 2},\\
\dsize\text{Re}\left( w_1^2{\partial^2u \over \partial w_1^2}\right)\ne 0,
\qquad & \text{when } |w_1|^2={1\over 2}, \\
\dsize\text{Im}\left( w_1^2{\partial^2u \over \partial w_1^2}\right)=0, \qquad
& \text{when } |w_1|^2 \le 1. 
\endcases 
\tag6.5
$$ 
In polar coordinates $w_1=re^{i\vartheta }$, (6.5) becomes
$$
\cases
\dsize (1-r^2H){\partial^2 u\over\partial r^2}+({1\over r}-2r+rH){\partial u
\over \partial r}+({1\over r^2}+H){\partial^2 u\over\partial \vartheta^2 }=0 ,
& \text{ when } r={\sqrt 2 \over 2},\\
\dsize r^2{\partial^2 u\over\partial r^2}-r{\partial u\over\partial r}-
{\partial^2u \over\partial\vartheta^2 }\ne 0, & \text{ when }r={\sqrt 2
\over 2}, \\
\dsize {\partial u\over \partial \vartheta}-r{\partial^2 u\over \partial
\vartheta \partial r}=0, & \text{ when } 0\le r \le 1.
\endcases
\tag6.6
$$
Any real function $u$ which is independent of $\vartheta$ and satisfies
$$
\cases
\dsize {\partial u\over \partial r}={H\over 2}-1, \\
\dsize {\partial^2u\over\partial r^2}={\sqrt 2\over 2}H,
\endcases
\qquad \text{when } r={\sqrt 2 \over 2}
\tag6.7
$$
solves the system (6.6). There are plenty of such real functions. For example,
$$
u=({H\over 2}-1)(r-{\sqrt 2\over 2})+{\sqrt 2\over 4}H
(r-{\sqrt 2\over 2})^2.
\tag6.8
$$
Therefore the proof is complete. \qed 
\enddemo 

\demo {Remark} No contact perturbation of $\text{id}_{S^3}: S_1^3 \to S_0^3$
with smooth Hamiltonian $u$ can reduce the magnitude of its Beltrami
tensor on $T_C$ at the level of the first variation. This fact becomes
clear if polar coordinates $w_1=re^{i\vartheta},\ \ w_2=\rho e^{i\varphi }$
are used to express
$$
\text{Im}(\overline\nu {\overline W_0}^2u)=2{\lambda^2-1 \over\lambda^2+1}
\left( -r{\partial^2u \over \partial r\partial \vartheta}-\rho {\partial^2u
\over \partial\rho\partial\vartheta}-r{\partial^2u\over\partial r\partial
\varphi}-\rho{\partial^2u\over\partial\rho\partial\varphi}+{\partial u\over
\partial\vartheta}+{\partial u\over\partial\varphi}\right). \tag 6.9
$$
In fact, the integral of right hand side of (6.9) over $(\vartheta,
\varphi)\in [0, 2\pi]\times [0, 2\pi]$ is zero for $u=u(\vartheta ,\varphi)$
is double $2\pi$-periodic in $(\vartheta ,\varphi)$. So $\text{Im}(\overline
\nu {\overline W_0}^2u)$ is neither positive nor negative on $T_C$. This
is the reason we need consider the second variation of $|\nu_{g_s}|$ to
demonstrate the symmetry breaking.
\enddemo

\enddocument